\documentclass[11pt,a4paper]{article}
\usepackage{amsmath}
\usepackage{enumerate}
\usepackage{amsfonts}
\usepackage{amsthm}
\usepackage{graphicx}
\usepackage[left=2.5cm,right=1.5cm,top=1cm,bottom=2cm]{geometry}

\usepackage{authblk}

\title{\textbf{Lie-B\"acklund symmetry  and non-invariant solutions of nonlinear evolution equations}}
\author[1]{I.M.Tsyfra}
\author[2]{W. Rzeszut}
\author[1]{V.A. Vladimirov}
\affil[1]{ AGH University of Science and Technology, Faculty of Applied Mathematics, 30 Mickiewicza Avenue, 30-059, Krakow, Poland}
\affil[2]{ Faculty of Mathematics and Computer Science, Jagiellonian University, Krakow, Poland}

\date{}

\begin{document}
\def\/*#1*/{}
\/*komentarz
afad\\
sdfsf\\
fgdfddgf
*/

\maketitle

\begin{abstract}
We study the symmetry reduction of nonlinear partial differential equations which are used for describing diffusion processes in nonhomogeneous medium. We find ansatzes reducing partial differential equations to systems of ordinary differential equations. The ansatzes are constructed by using operators of Lie-B\"acklund symmetry of third order ordinary differential equation. The method gives the possibility to find solutions which cannot be obtained by virtue of classical Lie method. Such solutions have been constructed for nonlinear diffusion equations which are invariant with respect to one-parameter, two-parameter and three-parameter Lie groups of point transformations.
\end{abstract}

\section{Introduction}

It is of common knowledge, that the most effective method for constructing solutions of nonlinear ODEs of mathematical physics is the symmetry reduction method which brings a PDE with several independent variables down to  another PDE with fewer independent variables, or even to ODE. The method can be both classical \cite{Olv} and non-classical \cite{BluCol,OlvRos1, OlvRos2,FusTsy}. In these cases the construction of a proper ansatz (by which we mean a general form of a invariant solution) boils down to solving a quasilinear first order DE, therefore ansatz includes one arbitrary function and the initial equation reduces to a single differential equation with fewer independent variables (especially an ODE).
In \cite{FokLiu} and \cite{Zhd} the concept of conditional  Lie-B\"acklund symmetry of evolution equations is proposed. By using this method one can reduce nonlinear evolution equations with two independent variables to system of ODEs.
The approach is used to construct exact solutions of nonlinear diffusion equations in \cite{Qu}. The relationship of generalized conditional symmetry of evolution equations to compatibility of system of differential equations is studied in \cite{KunPop}. Svirshchevskii \cite{Svirsh} put forward the reduction method for evolution equations of the form 
\[
u_t=K[u]
\]
where $u=u(t,x)$. The method is applicable if $K[u]\partial_u$ is a Lie-B\"acklund symmetry operator of a  linear homogeneous ODE. In this paper we use the method proposed in \cite{Tsy}, which is a generalization of  the Svirshchevskii's method, meaning we can analyze symmetries of a nonlinear (or nonhomogeneous linear) equation together with ODEs which include, besides dependent and independent variables, parametric variables and derivatives with respect to them. As it was shown in \cite{Tsy}, such generalization is important.
For example the equation
\[u_{xx}-a(t,x)u=0\]
is invariant w.r.t. the symmetry operator
\[\big(u_t-u_{xxx}+3\tfrac{u_{xx}u}{u_x}\big)\partial_u\]
provided that $a(t,x)$ is a solution to the KdV equation
\[a_t=-a_{xxx}+6aa_x.\]
It is clear then, that the method is related to the inverse scattering transformation method.\\
The idea is to use Lie-B\"acklund symmetries of ordinary differential equations (linear or nonlinear) for constructing solutions of evolutionary equations.
In this paper we consider a nonlinear evolutionary equation that describes transport phenomena in inhomogeneous  medium and apply a reduction method based on the symmetries  of  third order nonlinear ODEs.
Note that the proposed method could be applied to not only evolution-type equations, but to differential equations of any type as well. 

We present the results obtained for the model medium with exponential and polynomial heterogeneity. Within the method applied, nonlinear transport equation is reduced to a system of  three ODEs. After integrating (solving) the system of ODEs, we obtain exact solution of the initial equation.
Since the method applied differs from the classical Lie method, it does not enable   to construct algorithms for generation of new solutions, or production of conservation laws. Its only advantage is the preservation of the reduction property.  In addition, it doesn't ensure that none of the solutions obtained  could obtained within the classical method. Therefore there is a very important question of distinguishing truly new solutions obtained within the method proposed. 

Based on the fact that a set of point and Lie-B\"acklund symmetry operators (of the ODE) form a Lie algebra, we distinguish a class of diffusion equations  whose solutions, obtained with the help of the aforementioned approach,  cannot be obtained through the classical Lie method.

Furthermore, it can be used to construct a large class of nonlinear evolution equation all of which are reduced to systems of ordinary differential equations by the same ansatz and possess solutions which are not invariant in the classical Lie sense.

\section{Application of the method by using $3^{\text{rd}}$ order ODEs}
We consider third order differential equations of the following normal form
\[u_{xxx}-U(x,u,u_x,u_{xx})=0,\]
assuming that function $U$ is a power series of each of its arguments except $x$:
\[U=\sum_{j_0,j_1,j_2\in\mathbb{Z}}a_{j_0,j_1,j_2}(x)\,u^{j_0}\,u_x^{j_1}\,u_{xx}^{j_2}.\]
Here $a_{j_0,j_1,j_2}(x)$ are some smooth functions to be determined. We require that this ODE  admits  the symmetry operator $X=F(x,u,u_x,u_{xx})\partial_ u$, where $F$ is the right hand side of the transport equation, but the functions are independent of $t$, i.e., $u(x,t)\to\,u(x)$, $u_x(x,t)\to\,u_x(x)$, etc.
The corresponding symmetry condition reads as follows:
\begin{equation}
X^{(\infty)}\big(u_{xxx}-U(x,u,u_x,u_{xx})\big)\Big|_{u_{xxx}=U}=0,\label{eQ}\end{equation}
where $X^{(\infty)}$ is the prolongation of the Lie B\"{a}cklund symmetry generator $X$.
The algorithm of solving this equation is quite cumbersome and we omit it, presenting merely the summarizing results in Appendix.

\subsection{How to obtain  solutions of equation $u_t=(\tfrac{H(x)}{u})_{xx}+F(x,u,u_x)$?}

We use the fact that the PDE maintains the ansatz after we modify the PDE by adding the terms corresponding to symmetries of the ODE that generates that ansatz. Unfortunately, we are unable to present solutions to the majority of ODEs we deal within the method applied. Therefore we shall restrict ourselves to the ansatzes solvable in explicit forms, and connected with
$H(x)=\frac{1}{c_2x^2+c_1x+c_0}$. For the same reason, the choices of $F(x,u,u_x)$ are limited to those, for which the system of reduction equations is easily solvable in the real domain.
\\

For $H(x)=\frac{\kappa}{x}$, $(\kappa=const)$ equation
\begin{equation}
u_{xxx}=9\tfrac{u_{xx}u_{x}}{u}-12\tfrac{u_{x}^3}{u^2}
+\tfrac{6}{x}u_{xx}-\tfrac{18}{x}\tfrac{u_{x}^2}{u}-\tfrac{18}{x^2}u_{x}-\tfrac{12}{x^3}u
\label{gen1/x}
\end{equation}
admits the operator  $Q_1=(\tfrac{\kappa}{xu})_{xx}\partial_u$.
It has the solution
\[u(x)=\pm\tfrac{1}{x\sqrt{\varphi_2x^2+\varphi_1x+\varphi_0}}\]
and a $10$-parameter Lie group of point and contact symmetries:
\begin{eqnarray}
X_1=u \partial_u,\quad
X_2=xu_x \partial_u,\quad
X_3=x^2u^3 \partial_u,\quad
X_4=x^3u^3 \partial_u,\quad
X_5=x^4u^3 \partial_u,\quad
\end{eqnarray}
\begin{eqnarray}
X_6=(\tfrac{u}{x}+u_x) \partial_u,\quad
X_7=(2xu+x^2u_x) \partial_u,
\end{eqnarray}
\begin{eqnarray}
X_8=\tfrac{x^2u_x^2+4xuu_x+4u^2}{x^{2}u^{3}}\partial_u,\quad
X_9=\tfrac{x^2u_x^2+3xuu_x+2u^2}{x^{3}u^{3}} \partial_u,\quad
X_{10}=\tfrac{x^2u_x^2+2xuu_x+u^2}{x^{4}u^{3}} \partial_u.
\end{eqnarray}
\/*
The generalized symmetry of equation \eqref{gen1/x} has the form
$$x^2u^3\big(A_1(\zeta_0,\zeta_1,\zeta_2)+xA_2(\zeta_0,\zeta_1,\zeta_2)+x^2A_3(\zeta_0,\zeta_1,\zeta_2)\big)\frac{\partial}{\partial u},$$
where
$$\zeta_0=\tfrac{-x^2uu_{xx}+3x^2u_{x}^2+6xuu_x+6u^2}{x^2u^4},$$
$$\zeta_1=-2\tfrac{-x^2uu_{xx}+3x^2u_{x}^2+5xuu_x+4u^2}{x^3u^4},$$
$$\zeta_2=\tfrac{-x^2uu_{xx}+3x^2u_{x}^2+4xuu_x+3u^2}{x^4u^4}$$
are the invariants of equation \eqref{gen1/x} and $A_i$ $(i=1,2,3)$ are some smooth functions,
*/
\/*
so the reduction to a system of ODEs is possible for every equation from the class
\[u_t=\Big(\frac{\kappa}{xu}\Big)_{xx}+x^2u^3\big(A_1(\zeta_0,\zeta_1,\zeta_2)+xA_2(\zeta_0,\zeta_1,\zeta_2)+x^2A_3(\zeta_0,\zeta_1,\zeta_2)\big),\]
*/
so the reduction to a system of ODEs is possible for any equation fo the form
\[u_t=\Big(\frac{\kappa}{xu}\Big)_{xx}+a_i\sum_{i=1}^{10}X_iu,\quad a_i=const,~i=1,\ldots,10.\]
\/*
We want to modify the original equation by adding terms of order less than two, so we don't entirely change the type of the equation. We must add only characteristics of contact symmetries. For the purpose of this example we restrict ourselves to only few point symmetries.
*/
Let's consider equation
\[u_t=\Big(\frac{\kappa}{xu}\Big)_{xx}+a_1u+a_2xu_x+a_4(xu)^3,~\kappa,a_i\in\mathbb{R}.\]
Symmetry classification of this equation can be written down as
\[X_1=\partial_t,\]
\[X_2=2x\partial_x-3u\partial_u,\]
\[2a_1\neq 3a_2,~a_4=0 \Longrightarrow X_3=\mathrm{e}^{(2a_1-3a_2)t}\big(-a_2x\partial_x+\partial_t+a_1u\partial_u\big),\]
\[2a_1= 3a_2,~a_4=0 \Longrightarrow X_3=-a_2x\partial_x+t\partial_t+(a_1t+\tfrac{1}{2})u\partial_u.\]
Using the ansatz 
\[u(x,t)=\pm\tfrac{1}{x\sqrt{\varphi_2(t)x^2+\varphi_1(t)x+\varphi_0(t)}},\]
we get the following reduction equations
\[\varphi_2'+2(a_1-2a_2)\varphi_2=0,\]
\[\varphi_0'+2(a_1-a_2)\varphi_0=0,\]
\[\varphi_1'-\tfrac{1}{2}\kappa\varphi_1^2+2\kappa\varphi_0\varphi_2+2(a_1-\tfrac{3}{2}a_2)\varphi_1+2a_4=0.\]
For $2a_1\neq3a_2$, $a_4=0$
\[u(x,t)=\pm\frac{1}{x\sqrt{c_2\mathrm{e}^{2(2a_2-a_1)t}-2\sqrt{c_0c_2}\mathrm{e}^{(3a_2-2a_1)t}\tanh\big(\frac{\kappa\sqrt{c_0c_2}}{3a_2-2a_1}\mathrm{e}^{(3a_2-2a_1)t}\big) x+c_0\mathrm{e}^{2(a_2-a_1)t}}}.\]
For $a_1=\frac{3}{2}a_2=\frac{3}{2}a$, $a_4=0$
\[u(x,t)=\pm\frac{1}{x\sqrt{c_2\mathrm{e}^{at}x^2-2\sqrt{c_0c_2}\tanh\big(\kappa\sqrt{c_0c_2}(t+c_1)\big) x+c_0\mathrm{e}^{-at}}}.\]
For $a_1$, $a_2=0$ and $a_4$ arbitrary
\[u(x,t)=\pm\frac{1}{x\sqrt{c_2x^2-\frac{2\sqrt{\kappa^2c_0c_2+\kappa a_4}}{\kappa}\tanh\big(\sqrt{\kappa^2c_0c_2+\kappa a_4}(t+c_1)\big) x+c_0}}.\]
For $a_1\neq0$, $a_2\neq0$, $a_4\neq0$ we are unable to determine the solution.
Let us note, that none of the above solutions is invariant under any nonzero linear combination of $\{X_1,X_2,X_3\}$, as one can show that $$\sum_{i=1}^3\alpha_iX_i(u-u(x,t))\big|_{u=u(x,t)}=0 \Longrightarrow \forall_{j \in\{1,2,3\}}~\alpha_j=0.$$

let us consider the case $H(x)=\frac{\kappa}{x^2}$, $(\kappa=const).$ It can be verified by direct checking, that equation
\begin{equation}
u_{xxx}=9\tfrac{u_{xx}u_{x}}{u}-12\tfrac{u_{x}^3}{u^2}
+\tfrac{12}{x}u_{xx}-\tfrac{36}{x}\tfrac{u_{x}^2}{u}-\tfrac{60}{x^2}u_{x}-\tfrac{60}{x^3}u
\label{gen1/x^2}
\end{equation}
admits the operator  $Q_1=(\tfrac{\kappa}{x^2u})_{xx}\partial_u$.
It has a solution
\[u(x)=\pm\tfrac{1}{x^2\sqrt{\varphi_2x^2+\varphi_1x+\varphi_0}}\]
and a $10$-parameter Lie group of point and contact symmetries:
\begin{eqnarray}
X_1=u \partial_u,\quad
X_2=xu_x \partial_u,\quad
X_3=x^4u^3 \partial_u,\quad
X_4=x^5u^3 \partial_u,\quad
X_5=x^6u^3 \partial_u,\quad
\end{eqnarray}
\begin{eqnarray}
X_6=(2\tfrac{u}{x}+u_x) \partial_u,\quad
X_7=(3xu+x^2u_x) \partial_u,
\end{eqnarray}
\begin{eqnarray}
X_8=\tfrac{x^2u_x^2+6xuu_x+9u^2}{x^{4}u^{3}}\partial_u,\quad
X_9=\tfrac{x^2u_x^2+5xuu_x+6u^2}{x^{5}u^{3}} \partial_u,\quad
X_{10}=\tfrac{x^2u_x^2+4xuu_x+4u^2}{x^{6}u^{3}} \partial_u,
\end{eqnarray}
\/*
The generalized symmetry of equation \eqref{gen1/x^2} has the form
$$x^4u^3\big(B_1(\psi_0,\psi_1,\psi_2)+xB_2(\psi_0,\psi_1,\psi_2)+x^2B_3(\psi_0,\psi_1,\psi_2)\big)\frac{\partial}{\partial u},$$
where
$$\psi_0=\tfrac{-x^2uu_{xx}+3x^2u_{x}^2+10xuu_x+15u^2}{x^4u^4},$$
$$\psi_1=-2\tfrac{-x^2uu_{xx}+3x^2u_{x}^2+9xuu_x+12u^2}{x^5u^4},$$
$$\psi_2=\tfrac{-x^2uu_{xx}+3x^2u_{x}^2+8xuu_x+10u^2}{x^6u^4}$$
are the invariants of equation \eqref{gen1/x^2} and $B_i$ $(i=1,2,3)$ are some smooth functions, so the reduction to a system of ODEs is possible for every equation from the class
\[u_t=\Big(\frac{\kappa}{x^2u}\Big)_{xx}+x^2u^3\big(B_1(\psi_0,\psi_1,\psi_2)+xB_2(\psi_0,\psi_1,\psi_2)+x^2B_3(\psi_0,\psi_1,\psi_2)\big).\]
*/
 so the reduction to a system of ODEs is possible for every equation from the class
\[u_t=\Big(\frac{\kappa}{x^2u}\Big)_{xx}+a_i\sum_{i=1}^{10}X_iu,\quad a_i=const,~i=1,\ldots,10.\]
Let's consider equation
\[u_t=\Big(\frac{\kappa}{x^2u}\Big)_{xx}+a_4x^5u^3+a_5x^6u^3+a_7(3xu+x^2u_x),~\kappa,a_i\in\mathbb{R}.\]
Symmetry classification of this equation can be written down as
\[Y_1=\partial_t,\]
\[a_4\neq0\wedge a_7=0 \Longrightarrow Y_2=(-2x^2\tfrac{a_5}{a_4}-2x)\partial_x+t\partial_t+\tfrac{3}{2}(4x\tfrac{a_5}{a_4}+3)u\partial_u,\]
\[a_4=0 \Longrightarrow Y_2= -x\partial_x+t\partial_t+\tfrac{5}{2}u\partial_u,~Y_3=x^2\partial_x-3xu\partial_u,\]
\[a_4=a_5=0 \Longrightarrow Y_4=-a_7tx^2\partial_x+t\partial_t+(3a_7tx+\tfrac{1}{2})u\partial_u.\]
Since the ansatz is
\[u(x,t)=\pm\tfrac{1}{x^2\sqrt{\varphi_2(t)x^2+\varphi_1(t)x+\varphi_0(t)}},\]
the reduction equations are
\[\varphi_2'+2\kappa\varphi_0\varphi_2-\tfrac{1}{2}\kappa\varphi_1^2+2a_5+a_7\varphi_1=0,\]
\[\varphi_1'+2a_4+2a_7\varphi_0=0,\]
\[\varphi_0'=0,\]
with a solution
\[\varphi_2=\tfrac{(c_0a_7+a_4)^2}{c_0}t^2-(c_0a_7+a_4)(\tfrac{c_1}{c_0}+\tfrac{a_4}{\kappa c_0^2})t+\tfrac{c_1^2}{4c_0}+\tfrac{c_1a_4-2c_0a_5}{2\kappa c_0^2}+\tfrac{a_4(c_0a_7+a_4)}{2 \kappa^2 c_0^3}+c_2\mathrm{e}^{-2c_0\kappa t},\]
\[\varphi_1=-2(c_0a_7+a_4)t+c_1,\]
\[\varphi_0=c_0,\]
\[c_0\neq 0,\]
or
\[\varphi_2=\tfrac{2}{3}\kappa a_4^2t^3 +a_4(a_7-\kappa c_1)t^2+(\tfrac{1}{2}c_1^2\kappa-c_1a_7-2a_5)t+c_2,\]
\[\varphi_1=-2a_4t+c_1,\]
\[\varphi_0=0.\]
The  above solutions are invariant under the nontrivial linear combination of $\{Y_1,Y_2,Y_3,Y_4\}$ only when $a_4=a_5=0$ or $a_4=\varphi_0=0$.

Let's consider equation
\[u_t=\Big(\frac{\kappa}{x^2u}\Big)_{xx}+a_1u+a_4x^5u^3+a_5x^6u^3,~\kappa,a_i\in\mathbb{R},~a_1\neq0.\]
Symmetry classification of this equation can be written down as
\[Z_1=\partial_t,\]
\[a_4=0 \Longrightarrow Z_2=x^2\partial_x-3xu\partial_u,\]
\[a_4=a_5=0 \Longrightarrow Z_3=x\partial_x-2u\partial_u,~Z_4=\mathrm{e}^{2a_1t}(\partial_t+a_1u\partial_u).\]
Since the ansatz is
\[u(x,t)=\pm\tfrac{1}{x^2\sqrt{\varphi_2(t)x^2+\varphi_1(t)x+\varphi_0(t)}},\]
the reduction equations are as follows:
\[\varphi_2'+2\kappa\varphi_0\varphi_2-\tfrac{1}{2}\kappa\varphi_1^2+2a_1\varphi_2+2a_5=0,\]
\[\varphi_1'+2a_1\varphi_1+2a_4=0,\]
\[\varphi_0'+2a_1\varphi_0=0,\]
The bove system is fulfilled if
\[\varphi_2=\mathrm{e}^{-2a_1t}\bigg(\mathrm{e}^{\tfrac{c_0\kappa}{a_1}\mathrm{e}^{-2a_1t}}\Big(c_2-\tfrac{\kappa}{4a_1^2}\big(2c_1a_4-c_0(4a_5-\tfrac{a_4^2\kappa}{a_1^2})\big)~\Gamma\big(0,\tfrac{c_0\kappa}{a_1}\mathrm{e}^{-2a_1t}\big)\Big)+\tfrac{c_1^2}{4c_0}\bigg)-\tfrac{a_5}{a_1}+\tfrac{\kappa a_4^2}{4a_1^3},\]
\[\varphi_1=-\tfrac{a_4}{a_1}+c_1\mathrm{e}^{-2a_1t},\]
\[\varphi_0=c_0\mathrm{e}^{-2a_1t},\]
\[c_0\neq0,\]
or
\[\varphi_2=(c_2-\tfrac{c_1a_4\kappa}{a_1}t)\mathrm{e}^{-2a_1t}-\tfrac{\kappa c_1^2}{4a_1}\mathrm{e}^{-4a_1t}-\tfrac{a_5}{a_1}+\tfrac{\kappa a_4^2}{4a_1^3},\]
\[\varphi_1=-\tfrac{a_4}{a_1}+c_1\mathrm{e}^{-2a_1t},\]
\[\varphi_0=0.\]
The solution $u(x,t)$ with functions $\varphi_i$ as above is not invariant in classical sense
%under a nonzero linear combination of $\{Z_1,Z_2,Z_3,Z_4\}$ %
when $a_4\neq 0$ or $a_5\neq 0$.

Let's consider equation
\[u_t=\Big(\frac{\kappa}{x^2u}\Big)_{xx}+a_3x^4u^3+a_4x^5u^3+a_5x^6u^3,~\kappa,a_i\in\mathbb{R},~a_3\neq0.\]
Symmetry classification of this equation can be written down as
\[W_1=\partial_t,\]
\[a_4=a_5=0 \Longrightarrow W_2=x\partial_x-2u\partial_u.\]
Since the ansatz is
\[u(x,t)=\pm\tfrac{1}{x^2\sqrt{\varphi_2(t)x^2+\varphi_1(t)x+\varphi_0(t)}},\]
the the functions $\varphi_\nu,\,\,\,\nu=0,1,2$ satisfy the system
\[\varphi_2'+2\kappa\varphi_0\varphi_2-\tfrac{1}{2}\kappa\varphi_1^2+2a_5=0,\]
\[\varphi_1'+2a_4=0,\]
\[\varphi_0'+2a_3=0,\]
having the  solution
\[\varphi_2=\mathrm{e}^{2\kappa (a_3t^2-c_0 t)}\bigg(
\tfrac{-\sqrt{\pi}\big(\kappa(c_1-\tfrac{c_0a_4}{a_3})^2+(\tfrac{a_4^2}{a_3}-4a_5)\big)}{4\sqrt{2a_3\kappa}}\mathrm{e}^{\tfrac{\kappa c_0^2}{2a_3}}\text{erf}\big(\tfrac{\kappa(-2a_3t+c_0)}{\sqrt{2a_3\kappa}}\big)+c_2\bigg)+\tfrac{a_4}{4a_3}(-2a_4t+2c_1-\tfrac{c_0a_4}{a_3}),\]
\[\varphi_1=-2a_4t+c_1,\]
\[\varphi_0=-2a_3t+c_0.\]
The solution $u(x,t)$ with functions $\varphi_i$ as above is not invariant under the translations generated by $W_1=\partial_t$ when $a_4\neq$ or $a_5\neq0.$ It is not not invariant under  nonzero linear combination of $\{W_1,W_2\}$ for $a_4=a_5=0$.
%%%%%%%% to co poni�ej jest do wyrzucenia
\/*
\newline
To show that the property of reduction to a system of ODEs is not limited to the characteristics of point symmetries of an equation, we will now take contact symmetry operator $X_8$ and the equation
\[u_t=\Big(\frac{\kappa}{x^2u}\Big)_{xx}+a_8\tfrac{x^2u_x^2+6xuu_x+9u^2}{x^{4}u^{3}},~\kappa,a_8\in\mathbb{R},~a_8\neq0.\]
The reduction equations are
\[\varphi_0'+2a_8\varphi_0^2=0,\]
\[\varphi_1'+2a_8\varphi_0\varphi_1=0,\]
\[\varphi_2'+2\varphi_0\varphi_2=0,\]
with a solution
\[\varphi_0=\tfrac{1}{2a_8t+c_0},~\varphi_1=\tfrac{c_1}{2a_8t+c_0},~\varphi_2=\tfrac{\tfrac{1}{4}c_1^2}{2a_8t+c_0}+c_2(2a_8t+c_0)^{-\frac{1}{a_8}}\]
The equation has 4 point symmetry operators, meaning the solution
\[u(x,t)=\pm\frac{1}{x^2\sqrt{\tfrac{\big(1+\tfrac{1}{2}c_1x\big)^2}{2a_8t+c_0}+c_2(2a_8t+c_0)^{-\tfrac{1}{a_8}}}}\]
depending on three constants is invariant.
*/
\newpage
\section{Appendix}
\noindent Third order ODEs admitting Lie-B\"{a}cklund symmetries with infinitesimal generator $X=\big(\tfrac{H(x)}{u}\big)_{xx}\tfrac{\partial}{\partial u}$\\ ($k_1$, $k_2$, $k_3$, $k_4$, $k_5$, $\gamma$, $n$ are all real arbitrary constants, $H$ is differentiable) are as follows:\\
\noindent
\[u_{xxx}=7\tfrac{u_{xx}u_{x}}{u}-8\tfrac{u_{x}^3}{u^2}
-3\tfrac{H_{x}}{H}u_{xx}+8\tfrac{H_{x}}{H}\tfrac{u_{x}^2}{u}-4\tfrac{H_{xx}}{H}u_{x}+\tfrac{H_{xxx}}{H}u\]
for arbitrary $H$,
\[u_{xxx}=\tfrac{u_{xx}^2}{u_{x}}+4\tfrac{u_{xx}u_{x}}{u}-6\tfrac{u_{x}^3}{u^2}
-u_{xx}-\tfrac{uu_{xx}}{u_{x}}+6\tfrac{u_{x}^2}{u}-3u_{x}+u,
\]
for $H=\mathrm{e}^x$
\noindent
\[u_{xxx}=\tfrac{u_{xx}^2}{u_{x}}+4\tfrac{u_{xx}u_{x}}{u}-6\tfrac{u_{x}^3}{u^2}
+\tfrac{1-n}{x+\gamma}u_{xx}+\tfrac{n(1-n)}{(x+\gamma)^2}\tfrac{uu_{xx}}{u_{x}}+\tfrac{6n-2}{x+\gamma}\tfrac{u_{x}^2}{u}-\tfrac{n(3n-5)}{(x+\gamma)^2}u_{x}+\tfrac{n(n-1)(n-3)}{(x+\gamma)^3}u\]
for $H=(x+\gamma)^n$.
\noindent
Also if $H$ has the form $H=\tfrac{1}{c_2x^2+c_1x+c_0}$, then
\begin{center}
$
\begin{array}{rl}
u_{xxx}=&9\tfrac{u_{xx}u_{x}}{u}-12\tfrac{u_{x}^3}{u^2}
+\tfrac{6(2c_2x+c_1)}{c_2x^2+c_1x+c_0}u_{xx}
-\tfrac{18(2c_2x+c_1)}{c_2x^2+c_1x+c_0}\tfrac{u_{x}^2}{u}\\
&-\tfrac{6(10c_2^2x^2+10c_1c_2x-2c_0c_2+3c_1^2)}{(c_2x^2+c_1x+c_0)^2}u_{x}
-\tfrac{6(2c_2x+c_1)(5c_2^2x^2+5c_1c_2x-3c_0c_2+2c_1^2)}{(c_2x^2+c_1x+c_0)^3}u,
\end{array}
$
\end{center}
%\[u_{xxx}=9\tfrac{u_{xx}u_{x}}{u}-12\tfrac{u_{x}^3}{u^2}
%-6\tfrac{H_{x}}{H}u_{xx}+18\tfrac{H_{x}}{H}\tfrac{u_{x}^2}{u}-6\bigg(\tfrac{H_{xx}}{H}+\tfrac{(H_{x})^2}{H^2}\bigg)u_{x}+\bigg(\tfrac{1}{2}\tfrac{H_{xxx}}{H}+6\tfrac{H_{xx}H_{x}}{H^3}-3\tfrac{(H_{x})^3}{H^3}\bigg)u.\]
if $H=(x+\gamma)^{-\tfrac{1}{2}}$, then
\noindent
\[u_{xxx}=10\tfrac{u_{xx}u_{x}}{u}-15\tfrac{u_{x}^3}{u^2}
+ \tfrac{5k_2(x+\gamma)^2+3k_1}{(x+\gamma)\big(k_2(x+\gamma)^2+k_1\big)}u_{xx}
- \tfrac{35k_2(x+\gamma)^2+23k_1}{2(x+\gamma)\big(k_2(x+\gamma)^2+k_1\big)}\tfrac{u_{x}^2}{u}
\]\[
-\tfrac{25k_2(x+\gamma)^2+15k_1}{2(x+\gamma)^2\big(k_2(x+\gamma)^2+k_1\big)}u_{x}
-\tfrac{45k_2(x+\gamma)^2+25k_1}{8(x+\gamma)^3\big(k_2(x+\gamma)^2+k_1\big)}u+\tfrac{k_3(x+\gamma)^{3/2}}{k_2(x+\gamma)^2+k_1}u^4,\quad k_1\neq0 \vee k_2\neq0
\]
%\[u_{xxx}=10\tfrac{u_{xx}u_{x}}{u}-15\tfrac{u_{x}^3}{u^2}
%+ \tfrac{5}{x+\gamma}u_{xx}
%- \tfrac{35}{2(x+\gamma)}\tfrac{u_{x}^2}{u}
%-\tfrac{25}{2(x+\gamma)^2}u_{x}
%-\tfrac{45}{8(x+\gamma)^3}u+\tfrac{k_1}{(x+\gamma)^{1/2}}u^4,
%\]
%\[u_{xxx}=10\tfrac{u_{xx}u_{x}}{u}-15\tfrac{u_{x}^3}{u^2}
%+ \tfrac{3}{x+\gamma}u_{xx}
%- \tfrac{23}{2(x+\gamma)}\tfrac{u_{x}^2}{u}
%-\tfrac{15}{2(x+\gamma)^2}u_{x}
%-\tfrac{25}{8(x+\gamma)^3}u+k_1(x+\gamma)^{3/2}u^4,
%\]
\noindent
if $H=(x+\gamma)^{-\tfrac{3}{2}}$, then
\noindent
\[u_{xxx}=10\tfrac{u_{xx}u_{x}}{u}-15\tfrac{u_{x}^3}{u^2}
+ \tfrac{12k_2(x+\gamma)^2+10k_1}{(x+\gamma)\big(k_2(x+\gamma)^2+k_1\big)}u_{xx}
- \tfrac{87k_2(x+\gamma)^2+75k_1}{2(x+\gamma)\big(k_2(x+\gamma)^2+k_1\big)}\tfrac{u_{x}^2}{u}
\]\[
-\tfrac{135k_2(x+\gamma)^2+105k_1}{2(x+\gamma)^2\big(k_2(x+\gamma)^2+k_1\big)}u_{x}
-\tfrac{455k_2(x+\gamma)^2+315k_1}{8(x+\gamma)^3\big(k_2(x+\gamma)^2+k_1\big)}u+\tfrac{k_3(x+\gamma)^{7/2}}{k_2(x+\gamma)^2+k_1}u^4,\quad k_1\neq0\vee k_2\neq0
\]
%\[u_{xxx}=10\tfrac{u_{xx}u_{x}}{u}-15\tfrac{u_{x}^3}{u^2}
%+ \tfrac{12}{x+\gamma}u_{xx}
%- \tfrac{87}{2(x+\gamma)}\tfrac{u_{x}^2}{u}
%-\tfrac{135}{2(x+\gamma)^2}u_{x}
%-\tfrac{455}{8(x+\gamma)^3}u+k_1(x+\gamma)^{3/2}u^4,
%\]
%\[u_{xxx}=10\tfrac{u_{xx}u_{x}}{u}-15\tfrac{u_{x}^3}{u^2}
%+ \tfrac{10}{x+\gamma}u_{xx}
%- \tfrac{75}{2(x+\gamma)}\tfrac{u_{x}^2}{u}
%-\tfrac{105}{2(x+\gamma)^2}u_{x}
%-\tfrac{315}{8(x+\gamma)^3}u+k_1(x+\gamma)^{7/2}u^4,
%\]
\noindent
if $H=(x+\gamma)^{-2}$, then
\noindent
\begin{align*}
u_{xxx}=&9\tfrac{u_{xx}u_{x}}{u}-12\tfrac{u_{x}^3}{u^2}+\tfrac{12}{x+\gamma}u_{xx}-\tfrac{36}{x+\gamma}\tfrac{u_{x}^2}{u}-\tfrac{60}{(x+\gamma)^2}u_{x}-\tfrac{60}{(x+\gamma)^3}u\\
&+k_1((x+\gamma)^5u^3u_{x}+3(x+\gamma)^4u^4)+k_2((x+\gamma)^4u^3u_{x}+2(x+\gamma)^3u^4)\\
&+k_3((x+\gamma)uu_{xx}-2(x+\gamma)u_{x}^2-4uu_{x}-\tfrac{6u^2}{x+\gamma})
\end{align*}
\begin{align*}
u_{xxx}=&10\tfrac{u_{xx}u_{x}}{u}-15\tfrac{u_{x}^3}{u^2}
+ \tfrac{15}{x+\gamma}u_{xx}
- \tfrac{55}{x+\gamma}\tfrac{u_{x}^2}{u}
-\tfrac{105}{(x+\gamma)^2}u_{x}
-\tfrac{105}{(x+\gamma)^3}u\\
&+k_1\big((x+\gamma)uu_{xx}-3(x+\gamma)u_{x}^2-10uu_{x}-\tfrac{15}{x+\gamma}u^2\big)\\
&+k_2\big((x+\gamma)^2u^2u_{x}+3(x+\gamma)u^3\big)+k_3(x+\gamma)^3u^4+k_4(x+\gamma)^{5}u^5+k_5(x+\gamma)^6u^5,
\end{align*}
\begin{align*}
u_{xxx}=&10\tfrac{u_{xx}u_{x}}{u}-15\tfrac{u_{x}^3}{u^2}
+ \tfrac{14k_1+15k_2(x+\gamma)}{(x+\gamma)(k_1+k_2(x+\gamma))}u_{xx}
- \tfrac{52k_1+55k_2(x+\gamma)}{(x+\gamma)(k_1+k_2(x+\gamma))}\tfrac{u_{x}^2}{u}
-\tfrac{95k_1+105k_2(x+\gamma)}{(x+\gamma)^2(k_1+k_2(x+\gamma))}u_{x}\\
&-\tfrac{90k_1+105k_2(x+\gamma)}{(x+\gamma)^3(k_1+k_2(x+\gamma))}u+\tfrac{k_3(x+\gamma)^4}{k_1+k_2(x+\gamma)}u^4
+k_4\Big((x+\gamma)^2u^2u_{x}+\tfrac{2k_1+3k_2(x+\gamma)}{k_1+k_2(x+\gamma)}(x+\gamma)u^3\Big)\\
&+k_5\Big((x+\gamma)uu_{xx}-3(x+\gamma)u_{x}^2-\tfrac{9k_1+10k_2(x+\gamma)}{k_1+k_2(x+\gamma)}uu_{x}-\tfrac{12k_1+15k_2(x+\gamma)}{(x+\gamma)(k_1+k_2(x+\gamma))}u^2\Big)
,\quad k_1\neq0 \vee k_2\neq0,
%\begin{align*}
%u_{xxx}=&10\tfrac{u_{xx}u_{x}}{u}-15\tfrac{u_{x}^3}{u^2}
%+ \tfrac{14}{x+\gamma}u_{xx}
%- \tfrac{52}{x+\gamma}\tfrac{u_{x}^2}{u}
%-\tfrac{95}{(x+\gamma)^2}u_{x}
%-\tfrac{90}{(x+\gamma)^3}u+k_1(x+\gamma)^4u^4\\
%&+k_2\big((x+\gamma)^2u^2u_{x}+2(x+\gamma)u^3\big)
%+k_3\big((x+\gamma)uu_{xx}-3(x+\gamma)u_{x}^2-9uu_{x}-\tfrac{12u^2}{x+\gamma}\big)
\end{align*}
\noindent
\noindent
and finally, if $H=const$, then
\[u_{xxx}=10\tfrac{u_{xx}u_{x}}{u}-15\tfrac{u_{x}^3}{u^2}+\tfrac{u_{xx}}{x+k_1}-\tfrac{3}{x+k_1}\tfrac{u_{x}^2}{u}+k_2\big(u^2u_{x}+\tfrac{1}{x+k_1}u^3\big)+\tfrac{k_3}{x+k_1}u^4+k_4\big(uu_{xx}-3u_{x}^2-\tfrac{uu_{x}}{x+k_1}),\]
\[u_{xxx}=10\tfrac{u_{xx}u_{x}}{u}-15\tfrac{u_{x}^3}{u^2}+k_1u^2u_{x}+k_2u^5+k_3xu^5+k_4u^4+k_5\big(uu_{xx}-3u_{x}^2\big),\]
%\[u_{xxx}=\tfrac{u_{xx}^2}{u_{x}}+4\tfrac{u_{xx}u_{x}}{u}-6\tfrac{u_{x}^3}{u^2}+\tfrac{u_{xx}}{x+k_1}-\tfrac{2}{x+k_1}\tfrac{u_{x}^2}{u},\]
\[u_{xxx}=9\tfrac{u_{xx}u_{x}}{u}-12\tfrac{u_{x}^3}{u^2}+k_1(xu^3u_{x}+u^4)+k_2u^3u_{x},\]
\[u_{xxx}=\tfrac{u_{xx}^2}{u_{x}}+4\tfrac{u_{xx}u_{x}}{u}-6\tfrac{u_{x}^3}{u^2}.\]
%The ordinary equation $u_{xx}=9\frac{u_{xx}u_{x}}{u}-12\frac{u_{x}^3}{u^2}$ corresponds to a known higher-order summetry of the equation $u_t=(-u^{-1})_{xx}$ (see: Bluman, Kumei).
\noindent

\section{Conclusions}
We have constructed solutions of nonlinear evolution equations describing the diffusion processes in nonhomogeneous medium by using the generalization of Svirshchevskii method given in \cite{Tsy}.
We show that the method gives us the possibility to obtain solutions which are not invariant ones in the classical Lie sense. We use the Lie-B\"acklund symmetry operators of the third order ordinary differential equations. The corresponding ansatzes reduce nonlinear diffusion equations to the systems of three ordinary differential equations.
This way one is able to obtain the solutions which cannot be constructed by classical Lie method in the cases when the dimension of invariance Lie algebra is equal to 1,2,3. If the Lie algebra of Lie invariance group of the transport equation under consideration is four-dimensional, then the solutions obtained by using our method could also be obtained via classical Lie symmetry method as is shown in Section 2. These results agree with the ones given in \cite{Tsy1} if the solutions are found with the help of point conditional symmetry operators. The approach can be also applied to construct other classes of diffusion-type equations (and exact solutions) by using Lie-B\"acklund symmetry of other ordinary differential equations given in Appendix.

\end{document}